\documentclass{article}

\usepackage{amsmath,amssymb,amsthm}
\usepackage{tikz}
\usepackage{color}
\usepackage[toc]{appendix}
\usepackage{graphicx}
\usepackage{fancyhdr}
\usepackage{enumitem}
\usepackage{bbm}
\usepackage{parskip}
\usepackage{float}
\usepackage{chngpage}
\usepackage{calc}
\usepackage{bigints}
\usepackage{array}
\usepackage{booktabs}
\usepackage{rotating}
\usepackage{multirow}
\usepackage{adjustbox}
\usepackage{tabularx}
\usepackage{verbatim}
\usepackage{mathtools}
\usepackage{ragged2e}
\usepackage[makeroom]{cancel}
\usepackage{caption}
\usepackage{hyperref}
\usepackage{caption}
\usepackage{subcaption}
\usepackage{appendix}
\usepackage{pgfplots}
\pgfplotsset{compat=1.16}
\usepackage[english]{babel}
\usepackage{hyphenat}
\usepackage[makeindex]{imakeidx}

\usepackage{amsmath}

\usetikzlibrary{datavisualization}
\usetikzlibrary{matrix}
\usetikzlibrary{datavisualization.formats.functions}

\newtheorem{theorem}{Theorem}[section]

\newtheorem{example}[theorem]{Example}
\newtheorem{remark}[theorem]{Remark}
\newtheorem{assumption}[theorem]{Assumption}

\makeatletter
\def\section{\@startsection {section}{1}{\z@}{3.25ex plus 1ex minus
		.2ex}{1.5ex plus .2ex}{\large\bf}}
\def\subsection{\@startsection{subsection}{2}{\z@}{3.25ex plus 1ex minus
		.2ex}{1.5ex plus .2ex}{\normalsize\bf}}
\@addtoreset{equation}{section} 
\makeatother

\title{Maximum Likelihood With a Time Varying Parameter}

\author{Alberto Lanconelli\thanks{Dipartimento di Scienze Statistiche Paolo Fortunati, Università di Bologna, Bologna, Italy. \textbf{e-mail}: alberto.lanconelli2@unibo.it} \and Christopher S. A. Lauria\thanks{Dipartimento di Scienze Statistiche Paolo Fortunati, Università di Bologna, Bologna, Italy. \textbf{e-mail}: christopher.lauria2@unibo.it}}

\date{\today}

\begin{document}
	
\maketitle
	
\bigskip
	
\begin{abstract}
We consider the problem of tracking an unknown time varying parameter that characterizes the probabilistic evolution of a sequence of independent observations. To this aim, we propose a stochastic gradient descent-based recursive scheme in which the log-likelihood of the observations acts as time varying gain function. We prove convergence in mean-square error in a suitable neighbourhood of the unknown time varying parameter and illustrate the details of our findings in the case where data are generated from distributions belonging to the exponential family.   
\end{abstract}
	
	
	
\allowdisplaybreaks

\section{Introduction}	

When estimating unknown parameters in a dynamic model
the optimum solution to the parameter estimation problem may not remain constant. 
Specifically, the optimal values of the model parameters may change through time
because of the evolution of the underlying process: finding them is, in general, not straightforward. A survey of basic techniques for tracking the time-varying dynamics of a system is provided in \cite{LJUNG1997} where recursive algorithms in non-stationary stochastic optimization are analysed under different assumptions about the true system's variations, see also \cite{simonetto2020time} for a review in a purely deterministic setting. In \cite{Delyon} the problem of tracking the random drifting parameters of a linear regression system is tackled, and \cite{Zhu} builds a computable tracking error bound for how a stochastic approximation with constant gain keeps up with a non-stationary target. Successively, \cite{Wilson} introduces a framework for sequentially solving convex stochastic minimization problems, where the distance between successive minimizers is bounded. The minimization
problems are then solved by sequentially applying an optimization algorithm, such as stochastic gradient descent (SGD). In a similar setting, \cite{Cao} establishes an upper bound on the regret of a projected SGD algorithm with respect to the drift of the dynamic optima, while \cite{Cutler} provides novel non-asymptotic convergence guarantees for stochastic algorithms with iterate averaging. \\
We study time-varying stochastic optimization in a general statistical setting where we assume we are given a sequence of independent observations $\{ X_t \}_{t \in \mathbb{N}}$ with associated densities possessing a parameter that changes through time.
In such a framework a problem of interest concerns finding a useful estimator of the time varying parameter at a certain time $t$ - generalizing the classical problem of parameter estimation from the static setting to the time varying parameter setting.  
Ideally, one would like to find a sequence of estimators that track the time varying parameter through time as closely as possible. We show that, under some assumptions, utilizing the celebrated SGD algorithm \cite{Robbins} produces a sequence of estimators that will eventually track the time varying parameter - up to a neighborhood - as the number of observations increase. \\
Established in a general setting that intersects with the frameworks utilized in \cite{Cao}, \cite{Cutler} and \cite{Wilson}, our results differ from previous work mainly in one aspect: that our objective functions have the specific form of expected log likelihoods, a dissimilarity that will be exploited by utilizing their informational theoretical properties. \\
The work we present is also linked to the class of score driven models \cite{Creal}. Score driven models are a class of observation driven models (here we are using the terminology introduced by \cite{Cox1981}) that update the dynamics of the time varying parameter through the score of the conditional distribution of the observations.
Specifically, the same proof technique we utilize to obtain our result can be used to show that a -so called- Newton-score update \cite{Blasques}, with the parameter that multiplies the score appropriately chosen, will track the time varying parameter of interest trough time even under possible model misspecificaiton.  \\
A final way to interpret the results we present in this work is as robustness results for a one batch stochastic gradient procedure in the case we are incorrectly assuming that our observations are identically distributed. Indeed, the results show that even if we incorrectly assumed that the true parameter is static (we have IID observations) utilizing a stochastic gradient algorithm with a time dependent single sized batch to optimize the log-likelihood allows us to track the pseudo true time varying parameter up to a neighborhood if it is not moving wildly.\\
The paper is organised as follows: in Section 2 we list and discuss the assumptions of our framework and state the main result. We then present a class of examples given by the exponential family and discuss the performance of SGD with respect to the one observation maximum likelihood estimator at each time. In the third section we provide a detailed proof of our main result.
	
\section{Statement of the main result}\label{intro}

Let $\{ X_t \}_{t \in \mathbb{N}}$ be a sequence of independent $m$-dimensional random vectors defined on a common probability space $(\Omega,\mathcal{F},\mathbb{P})$. In the sequel we will write $\mathbb{E}[\cdot]$ for the expected value with respect to the probability measure $\mathbb{P}$, $\|\cdot\|$ for the Euclidean norm in $\mathbb{R}^d$ and $\|\cdot\|_{\mathbb{L}^2(\Omega)}$ for $\mathbb{E}[\|\cdot\|^2]^{\frac{1}{2}}$.\\
We assume that for any $t\in\mathbb{N}$ the random vector $X_t$ possesses a joint probability density function which depends on the $d$-dimensional parameter $\lambda_t^*$, in symbols $X_t\thicksim p(\cdot|\lambda_t^*)$. Our aim is to estimate the sequence $\{\lambda_t^*\}_{t\in\mathbb{N}}$ through the observed values $\{X_t\}_{t\in\mathbb{N}}$: To this aim we choose $\lambda_1\in\mathbb{R}^d$ and utilize the SGD algorithm
\begin{align} \label{sgd}
	\lambda_{t+1} := \lambda_t + \alpha \nabla_{\lambda} \ln p(X_{t} | \lambda_t ),\quad t\in\mathbb{N}.
\end{align}
Utilizing SGD to attempt to track $\lambda_t^* $ is motivated by the principle underlying classical maximum likelihood estimation: in fact, under some canonical assumptions we will present below, $\lambda_t^*$ will be the maximum of the expected log-likelihood $ \lambda \rightarrow \mathbb{E} \left[ \ln p( X_{t} | \lambda  )  \right]$. Thus, finding a sequence of estimators that track the time varying parameter as closely as possible is connected to finding the maxima of a sequence of expected log-likelihoods, a generalization of the classical static framework.
Since we have no direct access to the expected log-likelihoods, but only a singe observation for each time $t$, we categorize the problem as a time varying $\emph{stochastic}$ optimization problem. 

The assumptions we will require to obtain our result are the following.
\begin{assumption} [Smoothness of the log-likelihood] \label{ea}
	The function
	\begin{align}\label{function}
		\mathbb{R}^d\ni\lambda\mapsto \ln p(x|\lambda)
	\end{align}
	is twice continuously differentiable for all $x\in\mathbb{R}^m$; moreover,
	\begin{align*}
		\partial_{\lambda_i}\partial_{\lambda_j}\mathbb{E}\left[\ln p(X_t|\lambda)\right]=\mathbb{E}\left[\partial_{\lambda_j}\partial_{\lambda_j}\ln p(X_t|\lambda)\right],
	\end{align*}
	for all $i,j\in\{1,...,d\}$ and $t\in\mathbb{N}$.
\end{assumption} 

\begin{assumption}[Strong convexity]\label{sclg}
The function in \eqref{function} is strongly convex uniformly with respect to $x\in\mathbb{R}^m$: i.e., there exists a positive constant $\ell$ such that for all $x\in\mathbb{R}^m$ the matrix $\mathcal{H}_{\lambda}[-\ln p(x|\lambda)]-\ell I_d$ is positive semi-definite. Here, $\mathcal{H}_{\lambda}[-\ln p(x|\lambda)]$ stands for the Hessian matrix of the function in \eqref{function} while $I_d$ denotes the $d\times d$ identity matrix.
\end{assumption} 

\begin{assumption}[Lipschitz continuity of the gradient]  \label{sclg 2}
	The function
	\begin{align*}
		\mathbb{R}^d\ni\lambda\mapsto \nabla_{\lambda}\ln p(x|\lambda)
	\end{align*}
	is globally Lipschitz continuous uniformly with respect to $x\in\mathbb{R}^m$: i.e., there exists a positive constant $L$ such that for all $x\in\mathbb{R}^m$ we have
	\begin{align*}
		\|\nabla_{\lambda}\ln p(x|\xi_1)-\nabla_{\lambda}\ln p(x|\xi_2)\|\leq L \|\xi_1-\xi_2\|,\quad \xi_1,\xi_2\in\mathbb{R}^d.
	\end{align*}
\end{assumption} 

Assumptions \ref{sclg} and \ref{sclg 2} are classical in the optimization literature, see for instance \cite{Boyd} and\cite{Bottou2018}; we have utilized the versions of \cite{Nesterov}. We remark that Assumption \ref{sclg} may seem excessively restrictive at first glance, but we will present in Example \ref{ex} below a large family of examples where it holds.

\begin{remark} \label{leml}
Assumptions \ref{ea} and \ref{sclg 2} imply that
	\begin{align*}
		\mathtt{I}(\lambda_t^*)\leq dL,
	\end{align*}
	 where we have denoted $\mathtt{I}(\lambda_t^*):=\mathbb{E}[\|\nabla_{\lambda}\ln p(X_t|\lambda_t^*)\|^2]$, i.e. the trace of \emph{Fisher information matrix} of $X_t$.
In fact,
	\begin{align*}
		\mathtt{I}(\lambda_t^*)&=\mathbb{E}[\|\nabla_{\lambda}\ln p(X_t|\lambda_t^*)\|^2]=\sum_{j=1}^d\mathbb{E}[(\partial_{\lambda_j}\ln p(X_t|\lambda_t^*))^2]=-\sum_{j=1}^d\mathbb{E}[\partial^2_{\lambda_j}\ln p(X_t|\lambda_t^*)]\\
		&=\sum_{j=1}^d\mathbb{E}[\partial^2_{\lambda_j}(-\ln p(X_t|\lambda_t^*))]=\sum_{j=1}^d\mathbb{E}[\langle \mathcal{H}_{\lambda}(-\ln p(X_t|\lambda_t^*))e_j,e_j\rangle]\\
		&\leq\sum_{j=1}^d\mathbb{E}[\langle LI_de_j,e_j\rangle]=dL.
	\end{align*}

\end{remark}

We will use Remark \ref{leml} to bound the quantity $\mathbb{E}[\|\nabla_{\lambda}\ln p(X_t|\lambda_t)\|^2] $. In the general setting utilized in the optimization literature a bound on $\mathbb{E}[\|\nabla_{\lambda}\ln p(X_t|\lambda_t)\|^2] $ requires an extra assumption, see \cite{Bottou2018} and the discussion in \cite{Nguyen2018}. In our setting we manage to avoid this type of additional assumption thanks to the properties of the Fisher information matrix.\\
Our last assumption concerns the evolution of the time varying parameter $\{\lambda_t^*\}_{t\in\mathbb{N}}$.

\begin{assumption}[Lipschitz continuity of the true parameter] \label{ltvp}
There exists a positive constant $K$ such that
\begin{align*}
\| \lambda_{t+1}^* - \lambda_t^*\| \le K\quad\mbox{ for all $t\in\mathbb{N}$}.
\end{align*}
\end{assumption} 

Assumption \ref{ltvp} has been used throughout the literature, see for example \cite{simonetto2020time}, \cite{Cao} and \cite{Wilson}, since a limitation on the behavior of the sequence of true parameters values must be imposed to be able to track it. \\

We can now state our main theorem.

\begin{theorem} \label{t1} 
Let Assumptions \ref{ea}, \ref{sclg}, \ref{sclg 2} and \ref{ltvp} hold. Then, for $\alpha\in[\frac{1}{\ell+L}, \frac{1}{L}[$ running the SGD \eqref{sgd} we obtain
\begin{align}\label{limsup}
	\limsup_{t\to+\infty}\|\lambda_{t+1}-\lambda_t^*\|_{\mathbb{L}^2(\Omega)}\leq\frac{\varphi(\alpha,L)K+\alpha\sqrt{2dL}}{1-\varphi(\alpha,L)},
\end{align}
where $\varphi(\alpha,L):=\sqrt{1-2L\alpha+2L^2\alpha^2}$. Moreover, the minimum of the right hand side in \eqref{limsup} is attained at $\alpha=\frac{1}{\ell+L}$ and in this case the last inequality reads
\begin{align} \label{limsup 2}
	\limsup_{t\to+\infty}\|\lambda_{t+1}-\lambda_t^*\|_{\mathbb{L}^2(\Omega)}\leq\frac{K\sqrt{\ell^2 + L^2} + \sqrt{2dL}}{\ell + L - \sqrt{\ell^2 + L^2}}.
\end{align}
\end{theorem}

\begin{remark}
Notice that $\lambda_{t+1}$ depends on $X_1, X_2, \dots, X_t $, so as an estimator it is natural to compare it with $\lambda_t^* $.
\end{remark}

\begin{remark}
In the case of model misspecification, i.e. when the true distribution of the observations is not included in the parametric model $\{p(\cdot|\lambda)\}_{\lambda\in\mathbb{R}^d}$, the same proof technique can be utilized to show that the recursion \eqref{sgd} will track the so called \emph{pseudo-true} time varying parameter $\tilde{\lambda}_t$ which is defined as 
\begin{align*}
\tilde{\lambda}_t:=\arg\max_{\lambda\in\mathbb{R}^d}\mathbb{E}[\ln p(X_t|\lambda)].
\end{align*}
We recall that the pseudo-true time varying parameter $\tilde{\lambda}_t$ minimizes the Kullback Leiber divergence between the law of the data generating process and the model densities at each time $t$, see \cite{White} and \cite{akaike1973} for additional details. \\
The only technical difference in the proof is that Remark \ref{leml} can't be used since $\mathbb{E}[\|\nabla_{\lambda}\ln p(X_t|\tilde{\lambda}_t)\|^2]$ is no longer related to the Fisher information matrix of $X_t$. Thus, an additional assumption is needed to control $E[ \|\nabla_{\lambda}\ln p(X_{t}|\tilde{\lambda}_t)\|^2] $ but this is standard practice in the optimization literature, see \cite{Nguyen2018} for a discussion on this kind of assumption.   
\end{remark}

\begin{example} \label{ex}
	The exponential family in canonical form provides a class of natural examples where Theorem \ref{t1} holds. Take as the parameter of interest the natural parameter of a distribution belonging to the exponential family put in canonical form, i.e.
	\begin{align*}
		p(x|\lambda )=h(x)\exp\{\langle\lambda,T(x)\rangle-A(\lambda)\},\quad x\in\mathbb{R}^m
	\end{align*}
	where $h: \mathbb{R}^m \rightarrow \mathbb{R}$ is a non-negative function, $T: \mathbb{R}^m \rightarrow \mathbb{R}^d $ is a sufficient statistic and $ A: \mathbb{R}^d \rightarrow \mathbb{R}$ must be chosen so that  $p(x|\lambda)$ integrates to one. \\
	A standard result for exponential families, see for instance Theorem 1.6.3 in \cite{bickel2001}, is that $A$ is a convex function of $\lambda$; this fact together with identities
	\begin{align*}
		\nabla_{\lambda} \ln p(x|\lambda ) =  T(x) -\nabla_{\lambda} A( \lambda),
	\end{align*}
and
\begin{align*}
	\mathcal{H}_{\lambda}[-\ln p(x|\lambda )]=-\mathcal{H}_{\lambda} A(\lambda),
\end{align*}
implies that one can find, restricting if necessary the range of $\lambda$ (and hence of $\{\lambda_t^*\}_{t\in\mathbb{N}}$) to a suitable convex compact set $\Lambda$, the positive constants $l$ and $L$ required for the validity of Assumptions \ref{sclg}-\ref{sclg 2}.\\
Note that the restriction of the range of $\lambda$ to the convex compact set $\Lambda$ is carried out by simply modifying \eqref{sgd} as 
\begin{align*}
	\bar{\lambda}_{t+1} := \Pi_{\Lambda}\left( \bar{\lambda}_t +  \alpha \nabla_{\lambda} \ln p(X_{t} | \bar{\lambda}_t ) \right),\quad t\in\mathbb{N},
\end{align*}
where $\Pi_{\Lambda}$ denotes the orthogonal projection onto the set $\Lambda$. This alternative scheme
doesn't affect the validity of Theorem \ref{t1}; in fact, from the contraction property of $\Pi_{\Lambda}$ we get
\begin{align*} 
\| \bar{\lambda}_{t+1} - \lambda_t^*  \|^2 = \| \Pi_{\Lambda}( \bar{\lambda}_t +  \alpha \nabla_{\lambda} \ln p(X_{t} | \bar{\lambda}_t ) ) - \lambda_t^*  \|^2 \le \|  \bar{\lambda}_t +  \alpha \nabla_{\lambda} \ln p(X_{t} | \bar{\lambda}_t )  - \lambda_t^*  \|^2, 
\end{align*}
and this corresponds to the first step in the proof of Theorem \ref{t1} (see Section 3 below for more details).
\end{example}

An important question concerning applied settings is whether the estimator $\lambda_t$ defined in \eqref{sgd} performs asymptotically better than the maximum likelihood estimator $\hat{\lambda}_t$ calculated by optimizing the one observation log-likelihood $ \ln p(X_{t}|\lambda_{t})$. The following example will showcase that there are indeed cases when utilizing \eqref{sgd} is beneficial.

\begin{figure}[h]
	\centering
	\begin{subfigure}[c]{0.50\textwidth}
		\includegraphics[width=\linewidth]{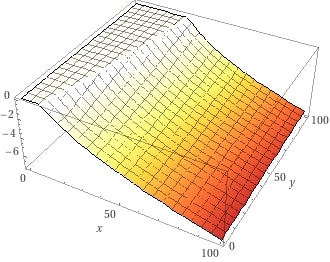}
		\caption{}
		\label{1}
	\end{subfigure}
	\begin{subfigure}[c]{0.40\textwidth}
		\includegraphics[width=\linewidth]{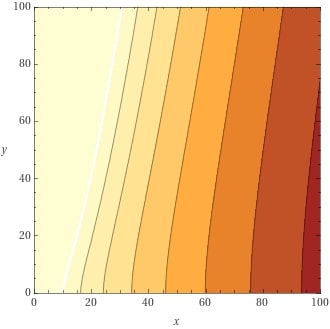}
		\caption{}
		\label{2}
	\end{subfigure}
\caption{Plot of the surface $z=\min\left\{\frac{K\sqrt{\ell^2 + L^2} + \sqrt{2L}}{\ell + L - \sqrt{\ell^2 + L^2}}-\sqrt{\ell},0\right\}$ from \eqref{ine} with $x=l$, $y=L-\ell$ and $K=1$.}
\label{1-2}

\end{figure}

\begin{figure}[H]
\begin{subfigure}[c]{0.50\textwidth}
		\includegraphics[width=\linewidth]{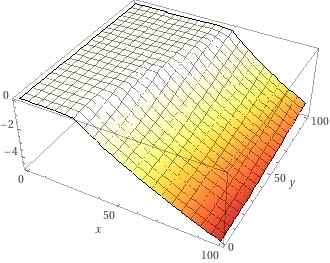}
		\caption{}
		\label{3}
	\end{subfigure}
	\begin{subfigure}[c]{0.40\textwidth}
		\includegraphics[width=\linewidth]{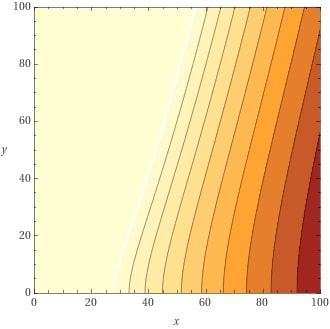}
		\caption{}
		\label{4}
	\end{subfigure}
\caption{Plot of the surface $z=\min\left\{\frac{K\sqrt{\ell^2 + L^2} + \sqrt{2L}}{\ell + L - \sqrt{\ell^2 + L^2}}-\sqrt{\ell},0\right\}$ from \eqref{ine} with $x=l$, $y=L-\ell$ and $K=2$.}
	\label{3-4}
	
\end{figure}

\begin{example}
Referring to Example \ref{ex} and setting $m=d=1$ for easiness of notation, we consider a sequence of independent observations $\{X_t\}_{t\in\mathbb{N}}$ with 
\begin{align*}
X_t\thicksim p(x|\lambda_{t}^* ):=h(x)\exp\{\lambda_t^* T(x)-A(\lambda_t^*)\},\quad x\in\mathbb{R}.
\end{align*}
We assume in addition that $\lambda\mapsto A''(\lambda)$ is continuous and we restrict the parameter space to $\Lambda = [ \lambda_m, \lambda_M  ] $ for suitable real numbers $\lambda_m< \lambda_M $. Observe that Assumptions \ref{sclg} and \ref{sclg 2} hold in this case with
\begin{align*}
	\ell = \min_{\lambda \in \Lambda} A''(\lambda), \quad L = \max_{\lambda \in \Lambda} A''(\lambda).
\end{align*}

In Theorem \ref{t1} we obtained an upper bound for the asymptotic mean-square error of $\lambda_t$ as defined in \eqref{sgd}. We now want to compare it with the mean-square error of the sufficient statistic $T(X_t)$, which we assume to be unbiased; this means considering the quantity
\begin{align}\label{w}
	\sqrt{\mathbb{E}[|T(X_t)-\lambda_t^*|^2]}=\sqrt{\mathbb{V}[T(X_t)]}=\sqrt{A''(\lambda_t^*)},
\end{align}
where the last equality follows from Theorem 1.6.2 in \cite{bickel2001}. Therefore, our estimator $\lambda_t$, performs asymptotically better than $T(X_t)$ if  
\begin{align}\label{z}
	 \frac{K\sqrt{\ell^2 + L^2} + \sqrt{2L}}{\ell + L - \sqrt{\ell^2 + L^2}} \le \sqrt{A''(\lambda_t^*)}\quad\mbox{ for all $t\in\mathbb{N}$}.
\end{align}
Here, the left hand side corresponds to right hand side in \eqref{limsup 2} with $d=1$ while the right hand side follows from \eqref{w}. We want this inequality to hold for all possible values of the sequence $\{\lambda_t^*\}_{t\in\mathbb{N}}$ and this is achieved by taking the infimum of the right hand side of \eqref{z}, i.e., we want
\begin{align}\label{ine} 
 \frac{K\sqrt{\ell^2 + L^2} + \sqrt{2L}}{\ell + L - \sqrt{\ell^2 + L^2}} \le \sqrt{\ell}.
\end{align}
A simple investigation of the previous inequality shows that the left hand side increases for small values of $\ell$ or large values of $L$; hence, there exist $\bar{\ell}$ and $\bar{L}$ such that for all $\bar{\ell}\leq \ell\leq L\leq \bar{L}$ the asymptotic mean-square error of $\lambda_t$ is lower than the mean-square error of the sufficient statistic $T(X_t)$. Figures \eqref{1-2} and\eqref{3-4} provide an illustration of this fact. Finally, notice that there are cases when the sufficient statistic of the exponential family is unbiased and coincides with the one observation maximum likelihood estimator, as is the case if we choose as the parameter of interest the variance of a Gaussian.

\end{example}

\newpage

\section{Proof of the main result}

Using \eqref{sgd} and expanding the squared Euclidian norm we can write  	
\begin{align}\label{a}
	\|\lambda_{t+1}-\lambda_t^*\|^2=&\|\lambda_t-\lambda_t^*+\alpha\nabla_{\lambda}\ln p(X_{t}|\lambda_t)\|^2\nonumber\\
		=&\|\lambda_t-\lambda_t^*\|^2+2\alpha\langle \lambda_t-\lambda_t^*,\nabla_{\lambda}\ln p(X_{t}|\lambda_t)\rangle+\alpha^2\|\nabla_{\lambda}\ln p(X_{t}|\lambda_t)\|^2\nonumber\\
		=&\|\lambda_t-\lambda_t^*\|^2+2\alpha\langle \lambda_t-\lambda_t^*,\nabla_{\lambda}\ln p(X_{t}|\lambda_t)-\nabla_{\lambda}\ln p(X_{t}|\lambda_t^*)\rangle\nonumber\\
		&+2\alpha\langle \lambda_t-\lambda_t^*,\nabla_{\lambda}\ln p(X_{t}|\lambda_t^*)\rangle+\alpha^2\|\nabla_{\lambda}\ln p(X_{t}|\lambda_t)\|^2\nonumber\\
		=&\|\lambda_t-\lambda_t^*\|^2+\mathcal{A}_1+2\alpha\langle \lambda_t-\lambda_t^*,\nabla_{\lambda}\ln p(X_{t}|\lambda_t^*)\rangle+\mathcal{A}_2,
\end{align}	
where we set
\begin{align*}
	\mathcal{A}_1:=2\alpha\langle \lambda_t-\lambda_t^*,\nabla_{\lambda}\ln p(X_{t}|\lambda_t)-\nabla_{\lambda}\ln p(X_{t}|\lambda_t^*)\rangle
\end{align*}
and
\begin{align*}
	\mathcal{A}_2:=\alpha^2\|\nabla_{\lambda}\ln p(X_{t}|\lambda_t)\|^2.
\end{align*}
To treat $\mathcal{A}_1$ we employ Theorem 2.1.12 from \cite{Nesterov}; with $C_1:=\frac{\ell L}{\ell+L}$ and $C_2=\frac{1}{\ell+L}$ this gives
\begin{align} \label{aa}
	\mathcal{A}_1\leq&-2\alpha C_1\|\lambda_t-\lambda_t^*\|^2-2\alpha C_2\|\nabla_{\lambda}\ln p(X_{t}|\lambda_t)-\nabla_{\lambda}\ln p(X_{t}|\lambda_t^*)\|^2;
\end{align}
moreover, using inequality $\|a+b\|^2\leq 2\|a\|^2+2\|b\|^2$ we get
\begin{align}\label{ac}
	\mathcal{A}_2\leq 2\alpha^2\|\nabla_{\lambda}\ln p(X_{t}|\lambda_t)-\nabla_{\lambda}\ln p(X_{t}|\lambda_t^*)\|^2+2\alpha^2\|\nabla_{\lambda}\ln p(X_{t}|\lambda_t^*)\|^2.
\end{align}
Combining \eqref{a} with \eqref{aa} and \eqref{ac} we obtain
\begin{align*}
	\|\lambda_{t+1}-\lambda_t^*\|^2\leq&(1-2\alpha C_1)\|\lambda_t-\lambda_t^*\|^2+2\alpha\langle \lambda_t-\lambda_t^*,\nabla_{\lambda}\ln p(X_{t}|\lambda_t^*)\rangle\\
	&+2\alpha(\alpha-C_2)\|\nabla_{\lambda}\ln p(X_{t}|\lambda_t)-\nabla_{\lambda}\ln p(X_{t}|\lambda_t^*)\|^2\\
	&+2\alpha^2\|\nabla_{\lambda}\ln p(X_{t}|\lambda_t^*)\|^2.
\end{align*}
Imposing that $2\alpha(\alpha-C_2)\geq 0$, or equivalently $\alpha\geq C_2$, we can utilize the Lipschitz continuity of the gradient in the second line above to get
\begin{align}\label{aab}
	\|\lambda_{t+1}-\lambda_t^*\|^2\leq&(1-2\alpha C_1+2\alpha(\alpha-C_2)L^2)\|\lambda_t-\lambda_t^*\|^2+2\alpha\langle \lambda_t-\lambda_t^*,\nabla_{\lambda}\ln p(X_{t}|\lambda_t^*)\rangle\nonumber\\
	&+2\alpha^2\|\nabla_{\lambda}\ln p(X_{t}|\lambda_t^*)\|^2.
\end{align}
Notice that according to the definitions of $C_1$ and $C_2$ we can write	
\begin{align*}
	1-2\alpha C_1+2L^2\alpha(\alpha-C_2)&=1-2\alpha(C_1+L^2C_2)+2L^2\alpha^2\\
	&=1-2\alpha\left(\frac{\ell L}{\ell+L}+\frac{L^2}{\ell+L}\right)+2L^2\alpha^2\\
	&=1-2L\alpha+2L^2\alpha^2.
\end{align*}
therefore, setting $\varphi(\alpha,L):=\sqrt{1-2L\alpha+2L^2\alpha^2}$ inequality \eqref{aab} now reads
\begin{align*}
	\|\lambda_{t+1}-\lambda_t^*\|^2\leq&\varphi(\alpha,L)^2\|\lambda_t-\lambda_t^*\|^2+2\alpha\langle \lambda_t-\lambda_t^*,\nabla_{\lambda}\ln p(X_{t}|\lambda_t^*)\rangle\nonumber\\
	&+2\alpha^2\|\nabla_{\lambda}\ln p(X_{t}|\lambda_t^*)\|^2.
\end{align*}
Taking the conditional expectation with respect to the sigma-algebra $\mathcal{F}_{t-1}:=\sigma(X_1,...,X_{t-1})$ of both sides above we obtain
\begin{align}\label{z}
	\mathbb{E}[\|\lambda_{t+1}-\lambda_t^*\|^2|\mathcal{F}_{t-1}]\leq&\varphi(\alpha,L)^2\|\lambda_t-\lambda_t^*\|^2+2\alpha\langle \lambda_t-\lambda_t^*,\mathbb{E}[\nabla_{\lambda}\ln p(X_{t}|\lambda_t^*)|\mathcal{F}_{t-1}]\rangle\nonumber\\
	&+2\alpha^2\mathbb{E}[\|\nabla_{\lambda}\ln p(X_{t}|\lambda_t^*)\|^2|\mathcal{F}_{t-1}]\nonumber\\
	=&\varphi(\alpha,L)^2\|\lambda_t-\lambda_t^*\|^2+2\alpha\langle \lambda_t-\lambda_t^*,\mathbb{E}[\nabla_{\lambda}\ln p(X_{t}|\lambda_t^*)]\rangle\nonumber\\
	&+2\alpha^2\mathbb{E}[\|\nabla_{\lambda}\ln p(X_{t}|\lambda_t^*)\|^2]\nonumber\\
	\leq &\varphi(\alpha,L)^2\|\lambda_t-\lambda_t^*\|^2+2\alpha^2dL.
\end{align}
Here, we have utilized that 
\begin{itemize}
	\item $\lambda_t$ is by construction $\mathcal{F}_{t-1}$-measurable for all $t\in\mathbb{N}$; 
	\item the $X_t$'s are independent;
	\item the expectation of the score is zero;
	\item Remark \ref{leml}.
\end{itemize}
We now compute the expectation of the first and last members of \eqref{z} to get
\begin{align*}
	\mathbb{E}[\|\lambda_{t+1}-\lambda_t^*\|^2]\leq\varphi(\alpha,L)^2\mathbb{E}[\|\lambda_t-\lambda_t^*\|^2]+2\alpha^2dL,
\end{align*}
which together with inequality $\sqrt{a+b}\leq \sqrt{a}+\sqrt{b}$ gives
\begin{align*}
	\|\lambda_{t+1}-\lambda_t^*\|_{\mathbb{L}^2(\Omega)}\leq\varphi(\alpha,L)\|\lambda_t-\lambda_t^*\|_{\mathbb{L}^2(\Omega)}+\alpha\sqrt{2dL}.
\end{align*}
The last step involves using Assumption \ref{ltvp} in the previous estimate to obtain
\begin{align*}
	\|\lambda_{t+1}-\lambda_t^*\|_{\mathbb{L}^2(\Omega)}\leq\varphi(\alpha,L)\|\lambda_t-\lambda_{t-1}^*\|_{\mathbb{L}^2(\Omega)}+\varphi(\alpha,L)K+\alpha\sqrt{2dL},
\end{align*}
which upon iteration yields
\begin{align*}
	\|\lambda_{t+1}-\lambda_t^*\|_{\mathbb{L}^2(\Omega)}\leq\varphi(\alpha,L)^{t-1}\|\lambda_2-\lambda_{1}^*\|_{\mathbb{L}^2(\Omega)}+(\varphi(\alpha,L)K+\alpha\sqrt{2dL})\frac{1-\varphi(\alpha,L)^{t-1}}{1-\varphi(\alpha,L)}.
\end{align*}
If $\alpha<\frac{1}{L}$, then $\varphi(\alpha,L)<1$; we can therefore take the limit as $t$ tends to infinity of both sides to get
\begin{align*}
	\limsup_{t\to+\infty}\|\lambda_{t+1}-\lambda_t^*\|_{\mathbb{L}^2(\Omega)}\leq\frac{\varphi(\alpha,L)K+\alpha\sqrt{2dL}}{1-\varphi(\alpha,L)};
\end{align*}
moreover, the minimum of the right hand side above is attained at $\alpha=\frac{1}{l+L}$ (in view of the constraints needed on $\alpha$ to recover inequality \eqref{aab}).

\bibliographystyle{apalike}
\bibliography{Time_varying_MLv6_arxiv}

\end{document}